\documentclass[10pt]{article}

\usepackage{amsmath,amsfonts,amsthm,amssymb,graphicx, longtable}
\usepackage[all,2cell,ps]{xy}

\theoremstyle{plain}
\newtheorem{thm}{Theorem}[section]

\newtheorem{prop}[thm]{Proposition}
\newtheorem{cor}[thm]{Corollary}

\theoremstyle{definition}

\newtheorem{ex}[thm]{Example}

\newtheorem*{rem}{Remark}

\theoremstyle{remark}

\newcommand{\bs}{\backslash}

\DeclareMathOperator{\U}{U}

\newcommand{\Z}{\mathbb Z}

\newcommand{\R}{\mathbb R}

\newcommand{\Q}{\mathbb Q}

\newcommand{\lam}{\lambda}
\newcommand{\Lam}{\Lambda}
\newcommand{\A}{\mathbb{A}}
\newcommand{\C}{\mathbb C}

\newcommand{\Gam}{\Gamma}

\newcommand{\conj}{\overline}

\newcommand{\Del}{\Delta}
\newcommand{\mc}{\mathcal}

\DeclareMathOperator{\SL}{SL}
\DeclareMathOperator{\PSL}{PSL}
\DeclareMathOperator{\GL}{GL}

\DeclareMathOperator{\PU}{PU}
\DeclareMathOperator{\SU}{SU}
\DeclareMathOperator{\PSU}{PSU}

\newcommand{\wt}{\widetilde}

\newenvironment{pf}{\begin{proof}}{\end{proof}}

\newenvironment{mat}{\left(\begin{matrix}}{\end{matrix}\right)}

\title{Volumes of Picard modular surfaces}

\author{Matthew Stover} 

\date{\today}

\begin{document}

\maketitle

\begin{abstract}

We show that the conjectural cusped complex hyperbolic 2-orbifolds of minimal volume are the two smallest arithmetic complex hyperbolic 2-orbifolds. We then show that every arithmetic cusped complex hyperbolic 2-manifold of minimal volume covers one of these two orbifolds. We also give all minimal volume manifolds that simultaneously cover both minimal orbifolds.

\end{abstract}



\section{Introduction}\label{intro}


Just as the modular group is a nonuniform arithmetic lattice in $\PSL_2(\R)$ acting on the hyperbolic plane, Picard modular groups are arithmetic lattices in $\PU(2,1)$ which act discretely and with finite covolume, though not cocompactly, on the complex hyperbolic plane. The corresponding algebraic surfaces are often called Picard modular surfaces. The purpose of this article is to study the minimal volume Picard modular surfaces, both in the orbifold and manifold categories.

Let $\Gam < \PU(2,1)$ be the fundamental group of a minimal Picard modular surface. Then $\Gam$ is a maximal lattice in $\PU(2,1)$ and is best described using Bruhat--Tits theory. This is done in $\S$\ref{arithmetic}. The preceding section, $\S$\ref{geometry}, is a brief review of complex hyperbolic space and the geometry of its finite volume quotients. In $\S$\ref{covolumes}, the noncompact arithmetic complex hyperbolic 2-orbifolds of minimal complex hyperbolic volume are determined using Prasad's formula \cite{Prasad}.


\begin{thm}\label{intro minimal orbifolds}

There are exactly two noncompact arithmetic complex hyperbolic $2$-orbifolds of complex hyperbolic volume $\pi^2 / 27$. They are commensurable and are the smallest volume Picard modular surfaces.

\end{thm}


We choose the metric for which the holomorphic sectional curvatures of the complex hyperbolic plane are $-1$. The proof uses computations done by Prasad--Yeung in their classification of fake projective planes \cite{Prasad--Yeung}. John Parker \cite{Parker} showed that $.25$ is a lower bound for the volume of a complex hyperbolic 2-orbifold and noted that these orbifolds come quite close to that bound. Conjecturally, these are the two smallest volume cusped complex hyperbolic 2-orbifolds.

A cusped complex hyperbolic manifold, i.e., a noncompact quotient of complex hyperbolic space with torsion-free fundamental group, must have Euler characteristic at least one. Therefore, Chern--Gauss--Bonnet implies that it has complex hyperbolic volume at least $8 \pi^2 / 3$. Parker found one such example \cite{Parker}. The last section, $\S$\ref{manifolds}, proves that his method of constructing such an example is, in a sense, the only way to do so.


\begin{thm}\label{intro manifold}

Let $M$ be an arithmetic complex hyperbolic $2$-manifold of volume $8 \pi^2 / 3$. Then $M$ covers one of the two orbifolds of Theorem \ref{intro minimal orbifolds}.

\end{thm}


The Appendix lists all the manifolds of Euler characteristic one that cover both orbifolds of Theorem \ref{intro minimal orbifolds}. This gives several new examples. For related work in hyperbolic 3-space, see \cite{CFJR}, and for hyperbolic 4-space, see \cite{Belolipetsky} (which considers all even-dimensional hyperbolic spaces) and \cite{Conder--Maclachlan}, the latter of which uses computational methods analogous to those used in the Appendix. In fact, the techniques of this paper are inspired by these three papers, along with \cite{Prasad--Yeung}. See the recent paper \cite{Belolipetsky--Emery} for odd-dimensional arithmetic hyperbolic orbifolds.


\section{Complex hyperbolic geometry}\label{geometry}


Let $V$ be a complex vector space of dimension three equipped with a Hermitian form $h$ of signature $(2,1)$. If $V_-$ is the set of $h$-negative vectors, set $\mathbf{H}_\C^2 = \mathbb{P}(V_-) \subset \mathbb{P}^2$. This is the \emph{complex hyperbolic plane}. Equipped with the metric associated to $h$, the group of biholomorphic isometries of $\mathbf{H}_\C^2$ is $\PU(2,1)$.

A \emph{complex hyperbolic} $2$-\emph{orbifold} is $\Gam \bs \mathbf{H}_\C^2$, where $\Gam < \PU(2,1)$ is a discrete subgroup of finite covolume. That is, the space $\Gam \bs \PU(2,1)$ has finite volume induced from the Haar measure on $\PU(2,1)$. The metric on $\mathbf{H}_\C^2$ also descends from the Haar measure on $\PU(2,1)$ under its identification with $\PU(2,1) / \U(2)$. Therefore, $\Gam$ has finite covolume in $\PU(2,1)$ if and only if $\Gam \bs \mathbf{H}_\C^2$ has finite complex hyperbolic volume.

Isometries of $\mathbf{H}_\C^2$ split into three categories based upon their action on $\mathbf{H}_\C^2$ and its ideal boundary. A nontrivial isometry is called \emph{elliptic} if it has a fixed point in the complex hyperbolic plane; it is \emph{loxodromic} if it is not elliptic and has two fixed points on the ideal boundary; it is \emph{parabolic} if it is not elliptic and fixes exactly one point on the ideal boundary.

The space $\Gam \bs \mathbf{H}_\C^2$ is a manifold if and only if $\Gam$ is torsion-free, which holds if and only if $\Gam$ contains no elliptic elements. That every complex hyperbolic 2-orbifold has a finite-sheeted manifold cover follows from Selberg's Lemma (see \cite{Alperin} for a proof). A lattice $\Gam$ is called \emph{uniform} if $\Gam \bs \mathbf{H}_\C^2$, or equivalently $\Gam \bs \PU(2,1)$, is compact and \emph{nonuniform} otherwise. A lattice is nonuniform if and only if it contains parabolic elements.

See \cite{Eberlein} for generalities on lattices and their corresponding quotient spaces. For the construction of complex hyperbolic space and its isometries, see \cite{Goldman}. The best-known examples of noncompact complex hyperbolic 2-orbifolds are \emph{Picard modular surfaces}. The corresponding lattices are precisely the nonuniform arithmetic lattices in $\PU(2,1)$. The next section describes these lattices in detail.


\section{Picard modular groups}\label{arithmetic}



\subsection{}\label{algebraic groups}


Let $d$ be a square-free positive integer and $k$ the imaginary quadratic field $\Q(\sqrt{-d})$. If $h$ is a Hermitian form on $k^3$ of signature $(2,1)$, the special unitary group of $h$ is a $\Q$-form of $\SU(2,1)$. That is, if $\mc{G}_h$ is the $\Q$-algebraic group such that
\[
\mc{G}_h(\Q) \cong \{ A \in \SL_3(k)\ :\ {}^t \conj{A} h A = h \},
\]
then $\mc{G}_h(\R) \cong \SU(2,1)$. Hermitian forms on $k^3$ of signature $(2,1)$ with the same determinant are isomorphic \cite[Ch.\ 10]{Scharlau}, from which it follows that all the $\mc{G}_h$ over a fixed $k$ are $\Q$-isomorphic \cite[$\S$1.2]{Prasad--Yeung}. In their notation, $k$ is $\Q$, $\ell$ is our imaginary quadratic field, and their division algebra $\mc{D}$ is $\ell$. They prove that we may assume that $h$ has determinant one. Since there is exactly one isomorphism class of Hermitian forms of determinant one on $\mc{D}^3$, this proves that each imaginary quadratic field determines one, and only one, $\Q$-form of $\SU(2,1)$.

The subgroup
\[
\mc{G}_h(\Z) = \{ A \in \SL_3(\mc{O}_k)\ :\ {}^t \conj{A} h A = h \},
\]
where $\mc{O}_k$ is the ring of integers of $k$, projects to a nonuniform lattice in $\PU(2,1)$. If $h$ are $h^\prime$ are two Hermitian forms over imaginary quadratic fields $k$ and $k^\prime$, then $\mc{G}_h(\Z)$ is commensurable with $\mc{G}_{h^\prime}(\Z)$ if and only if $k = k^\prime$, though the two groups are not necessarily isomorphic. In other words, to each imaginary quadratic field $k$, there corresponds a unique commensurability class of nonuniform arithmetic lattices in $\PU(2,1)$. This construction determines all commensurability classes of nonuniform arithmetic lattices in $\SU(2,1)$.


\subsection{}\label{arithmetic lattices}


Fix an imaginary quadratic field $k$ with integer ring $\mc{O}_k$, let $h$ be a Hermitian form on $k^3$, and let $\mc{G}$ be the associated $\Q$-form of $\SU(2,1)$. For each rational prime $p$, consider the group $\mc{G}(\Q_p)$ of $\Q_p$-points of $\mc{G}$. A direct computation shows that $\mc{G}(\Q_p)$ is isomorphic to $\SL_3(\Q_p)$ when $p$ splits in $\mc{O}_k$. If $p$ is inert or ramified, let $k_p$ denote the completion of $k$ at the prime ideal of $\mc{O}_k$ above $p$. This is a quadratic extension of $\Q_p$, and $\mc{G}(\Q_p)$ is the unitary group of $h$ for this extension.

For each $p$, let $K_p$ be a compact open subgroup of $\mc{G}(\Q_p)$, and suppose that $K_p$ is hyperspecial for all but finitely many $p$ (see \cite[$\S$3.8]{Tits}). This ensures that
\[
K_f = \prod_{p} K_p < \mc{G}(\mathbb{A}_f)
\]
is a compact open subgroup, where $\mathbb{A}_f$ is the finite adele ring of $\Q$. Embedding $\mc{G}(\Q)$ in $\mc{G}(\mathbb{A}_f)$ diagonally, $\Gam_{K_f} = K_f \cap \mc{G}(\Q)$ is a lattice in $\SU(2,1)$ commensurable with $\mc{G}(\Z)$. In fact, $\mc{G}(\Z)$ is amongst these lattices, where $K_p = \mc{G}(\Z_p)$ for all $p$. The integral structure is given by a choice of Hermitian form $h$ for some basis for $k^3$, as in $\S$\ref{algebraic groups}, and the group $\mc{G}(\Z_p)$ is hyperspecial for all $p$ under our choice of $h$. We denote $\mc{G}(\Z)$ by $\Gam_{\textrm{std}}$. This is the standard Picard modular group.

The conjugacy classes of maximal compact open subgroups of $\mc{G}(\Q_p)$ are determined by the vertices of its local Dynkin diagram. The cases of interest to us are worked out in detail in the examples of \cite{Tits}.

When $p$ splits in $\mc{O}_k$, let $K_p^{v_j}$, $j = 0, 1, 2$, be representatives for the conjugacy classes of maximal compact open subgroups, i.e., the maximal parahoric subgroups, of $\mc{G}(\Q_p)$. These are the stabilizers of three vertices, $v_0$, $v_1$, and $v_2$, of the corresponding Bruhat-Tits building. These vertices are $\GL_3(\Q_p)$-conjugate, but not $\SL_3(\Q_p)$-conjugate. They determine a chamber $t_0$ of the building, which is a triangle. Let $K_p^0$ be the associated Iwahori subgroup $K_p^{v_0} \cap K_p^{v_1} \cap K_p^{v_2}$.

If $p$ is inert or ramified, the building is a tree. There are two conjugacy classes of compact open subgroups, $K_p^{v_0}$  and $K_p^{v_1}$. We assume that $v_0$ corresponds to the hyperspecial vertex when $p$ is inert. The vertices $v_0$ and $v_1$ are adjoined by an edge $e_0$, which is a chamber of the building with associated stabilizer the Iwahori subgroup $K_p^0$. When $p$ is inert, the tree is not regular, and the groups $K_p^{v_j}$ are not conjugate. When $p$ is ramified, the groups are still not equivalent, but it is more subtle (see \cite[$\S$3.5]{Tits} and \cite[$\S$2.2]{Prasad--Yeung}).


\subsection{}\label{principal arithmetic lattices}


Since $\mc{G}$ is a $\Q$-form of $\SU(2,1)$, it has Strong Approximation \cite[Theorem 7.12]{Platonov--Rapinchuk}. Therefore, for any compact open $K_f < \mc{G}(\mathbb{A}_f)$, there exists $g \in \mc{G}(\Q)$ so that each factor of $g K_f g^{-1}$ is contained in one of the $K_p^{v_j}$ of $\S$\ref{arithmetic lattices}. Furthermore, the action of $\mc{G}(\Q_p)$ on the Bruhat--Tits building at $p$ is \emph{special} \cite[$\S$3.2]{Tits}, from which it follows that if $\Gam < \mc{G}(\Q)$ is a maximal lattice, then $\Gam = \Gam_{K_f}$ for some compact open subgroup of $\mc{G}(\A_f)$, and each factor of $K_f$ stabilizes a vertex of the building. Therefore, we can, and do, assume that any $\Gam_{K_f} < \mc{G}(\Q)$ is such that each factor is contained in one of the groups $K_p^{v_j}$.

However, the $\Gam_{K_f}$ which are maximal in $\mc{G}(\Q)$ do not necessarily determine the maximal lattices in $\mc{G}(\R) \cong \SU(2,1)$. Let $\conj{\mc{G}}$ be the corresponding $\Q$-form of $\PU(2,1)$. Then every maximal arithmetic lattice in $\PU(2,1)$ is contained in $\conj{\mc{G}}(\Q)$ \cite[Proposition 4.2]{Platonov--Rapinchuk}, but the projection from $\mc{G}(\Q)$ to $\conj{\mc{G}}(\Q)$ is not necessarily surjective. One can describe this phenomenon explicitly via lifts to the general unitary group. The lift $g$ of some element of $\conj{\mc{G}}(\Q)$ to $\mathrm{GU}(2,1)$ will have a special unitary representative lying in $\mc{G}(\Q)$ if and only if $\det(g)$ is a cube in $k$.

However, each maximal lattice in $\SU(2,1)$ in our commensurability class is the normalizer in $\mc{G}(\conj{\Q})$ of some $\Gam_{K_f}$ \cite[Proposition 1.4]{Borel--Prasad}. We may also assume that each factor of $K_f$ is either one of the vertex stabilizers $K_p^{v_j}$, or is an Iwahori subgroup $K_p^0$. Even further, when $p$ is inert or ramified, every element of $\mc{G}(\conj{\Q}_p)$ acts by special automorphisms, so we may assume that $K_f$ is Iwahori only at split primes.

Let $\conj{\Gam} < \mc{G}(\conj{\Q})$ be a maximal lattice, normalizing the lattice $\Gam_{K_f}$. Let $\mc{I}$ be the set of (split) primes for which $K_f$ is Iwahori at $p$. Then $K_f$ is one of the groups $K_p^{v_j}$ for all other $p$. Then, we have the following upper bound for $[\conj{\Gam} : \Gam_{K_f}]$, which will be of use several times throughout this paper.


\begin{prop}[See $\S$5.3 in \cite{Borel--Prasad} and Equation (0) in \cite{Prasad--Yeung}]\label{normalizer index}

Let $\conj{\Gam}$ and $\Gam_{K_f}$ be as above. Then
\[
[\conj{\Gam} : \Gam_{K_f}] \leq 3^{1 + |\mc{I}|} h_{k, 3},
\]
where $h_{k, 3}$ is the order of the $3$-primary part of the class group of the field $k$ from which the algebraic group $\mc{G}$ is defined.

\end{prop}



\section{Covolumes of Picard modular groups}\label{covolumes}



\subsection{}\label{volume formula}


Let $K_f$ and $\Gam_{K_f}$ be as in $\S$\ref{arithmetic}, that is, conjugate $K_f$ so that each factor is contained in one of the compact open subgroups from there. Prasad \cite{Prasad} gave an explicit formula for the covolume of arithmetic lattices in semisimple groups, which in our case determines the volume of $\Gam_{K_f} \bs \SU(2,1)$ with respect to the Haar measure $\mu$ on $\SU(2,1)$. Normalize $\mu$ such that the compact dual, $\mathbb{P}^2$, of $\SU(2,1)$ has Euler-Poincar\'{e} characteristic three and let $\mathrm{L}_k$ denote the $\mathrm{L}$-function of the quadratic extension $k/\Q$. We note that this normalization is used in \cite{Prasad--Yeung}, but differs from the normalization of \cite{Prasad}.


\begin{thm}\label{prasad's formula}

Let $K_f = \prod K_p$ be as in $\S$\ref{arithmetic} and $\Gam_{K_f} < \SU(2,1)$ the associated principal arithmetic lattice. Then
\[
\mu(\Gam_{K_f} \bs \SU(2,1)) = -\frac{1}{48} \mathrm{L}_k(-2) \prod_{p \in \mc{P}} \lam_p,
\]
where the $\lam_p$ are
\begin{align}
(p^2 + p + 1) (p + 1),\quad & p\ \textrm{split and}\ K_p = K_p^0, \nonumber \\
p^3 + 1,\quad & p\ \textrm{inert and}\ K_p = K_p^0, \nonumber \\
p + 1,\quad & p\ \textrm{ramified and}\ K_p = K_p^0, \nonumber \\
p^2 - p + 1,\quad & p\ \textrm{inert and}\ K_p = K_p^{v_1}, \nonumber
\end{align}
and $\lam_p=1$ otherwise.

\end{thm}


See \cite{Prasad--Yeung} for an explanation. The covolumes of the standard Picard modular groups, along with those of some congruence subgroups, were also computed by Holzapfel \cite{Holzapfel} and in the dissertation of Zeltinger \cite{Zeltinger}. By Proposition \ref{normalizer index}, it suffices to assume that $K_f$ is of the form $K_p^{v_j}$, i.e., is a maximal parahoric subgroup of $\mc{G}(\Q_p)$, when $p$ is inert or ramified, and is maximal parahoric at all split $p$, except for a finite set $\mc{I}$ at which it is the Iwahori subgroup $K_p^0$.

For every imaginary quadratic field $k$, the lattice $\Gam_{\mathrm{std}} = \mc{G}(\Z)$ has minimal covolume amongst principal arithmetic lattices in its commensurability class, since $\lam_p = 1$ for all $p$. However, $\Gam_{\mathrm{std}}$ is not the unique principal arithmetic lattice of minimal covolume in its commensurability class. It has `sister' lattices. The choice of a maximal compact subgroup at the ramified primes, i.e., the choice of $K_p^{v_0}$ versus $K_p^{v_1}$, has no effect on Prasad's formula, so there are at least $2^{|\mathrm{Ram}(k)|}$ distinct lattices of the same covolume, where $\mathrm{Ram}(k)$ is the set of rational primes which ramify in $k$. When the class number of $k$ is one, there are exactly $2^{|\mathrm{Ram}(k)|}$ such lattices \cite[$\S$5.2]{Prasad--Yeung}. A prime ramifies if and only if it divides the discriminant $d_k$ of $k$. Since $|d_k| > 1$, sisters always exist.


\subsection{}\label{upper bound}


We begin by examining the consequences of Proposition \ref{normalizer index} for a fixed commensurability class. That is, fix $d$ and let $\conj{\Gam}$ be a maximal arithmetic subgroup of $\mc{G}(\conj{\Q})$. Let $\Gam_{K_f} < \mc{G}(\Q)$ be the principal arithmetic lattice it normalizes. Then
\[
\mu(\conj{\Gam} \bs \SU(2,1)) = \frac{1}{[\conj{\Gam} : \Gam_{K_f}]} \mu(\Gam_{K_f} \bs \SU(2,1)) \geq -\frac{\prod_{p \in \mc{P}} \lam_p}{3^{1 + |\mc{I}|} h_{k,3}} \frac{\mathrm{L}_k(-2)}{48}.
\]

Notice that $\lam_p / 9 > 1$ for all $p \in \mc{I}$ and $\lam_p / 3 > 1$ for all $p > 2$ for which $K_p$ is Iwahori, regardless of the decomposition of $p$ in $k$. The same holds when $p$ is inert and $K_p$ is $K_p^{v_1}$. Assume $d \neq 3$. Then
\[
\mu(\conj{\Gam} \bs \SU(2,1)) \geq -\frac{1}{144 h_{k,3}} \mathrm{L}_k(-2) \geq -\frac{1}{144 h_k} \mathrm{L}_k(-2).
\]
Now, apply the functional equation
\[
\mathrm{L}_k(-2) = -\frac{|d_k|^{\frac{5}{2}}}{2 \pi^3} \mathrm{L}_k(3)
\]
to obtain
\[
\mu(\conj{\Gam} \bs \SU(2,1)) \geq \frac{|d_k|^{\frac{5}{2}}}{288 \pi^3 h_k} \mathrm{L}_k(3),
\]
where $d_k$ is the discriminant of $k$.

The proof of the Brauer--Siegel Theorem implies that
\[
h_k \leq w_d \frac{n (n - 1) (n - 1)!}{(2 \pi)^n} |d_k|^{\frac{n}{2}} \zeta(n) \mathrm{L}_k(n)
\]
for any $n > 1$, where $w_d$ is the number of roots of unity in $k$ and $\zeta(s)$ is the Riemann zeta function. We briefly explain how to derive this bound from \cite{Lang}. This inequality is a restatement of the inequality in the proof of Lemma 1 in $\S$1 of Chapter XVI. The term $\kappa$ is defined in Proposition 9 in $\S$7 of Chapter XIV. We recall that the regulator, $R$, is one for imaginary quadratic fields; $r_1$, the number of real places is one; and $r_2$, the number of complex conjugate places, is also one. Also, Lang uses the Dedekind zeta function $\zeta_k(s)$, but $\zeta_k(s) = \zeta(s) L_k(s)$.

If $d \neq 3,4$, then $w_d = 2$. In particular, for $n = 3$ and $d \neq 3,4$,
\[
\mu(\conj{\Gam} \bs \SU(2,1)) \geq \frac{|d_k|}{864 \zeta(3)}.
\]
Since $\zeta(3) < 5 / 4$, we get the following.


\begin{thm}\label{volume upper bound}

Let $k = \Q(\sqrt{-d})$ for $d \neq 1, 3$ and let $\conj{\Gam} < \SU(2,1)$ be a maximal Picard modular group defined via the imaginary quadratic field $k$. Then
\[
\mu(\conj{\Gam} \bs \SU(2,1)) \geq \frac{|d_k|}{1080},
\]
where $d_k$ is the discriminant of $k$.

\end{thm}


The bounds for $d = 1, 3$ are $1 / 32$ and $1 / 216$, respectively. The Chern--Gauss--Bonnet Theorem gives the following corollary (see $\S$\ref{orbifolds} and \cite{Prasad--Yeung}).


\begin{cor}\label{bound corollary}

Suppose $k = \Q(\sqrt{-d})$ for $d \neq 1, 3$ and that $M$ is a Picard modular surface with associated imaginary quadratic field $k$. Then the complex hyperbolic volume of $M$ is at least $|d_k| \pi^2 / 405$, where $d_k$ is the discriminant of $k$.

\end{cor}



\subsection{}\label{minimal}


For this section, let $\Gam_d = \Gam_{\mathrm{std}, d}$ be the standard Picard modular group for $\Q(\sqrt{-d})$.


\begin{thm}\label{3 is minimal}

Let $\conj{\Gam}_3$ be the normalizer of $\Gam_3$ in $\SU(2,1)$. Then $\conj{\Gam}_3$ is a three-fold extension of $\Gam_3$. It, and its sister, are the nonuniform arithmetic lattices in $\SU(2,1)$ of minimal covolume.

\end{thm}



\begin{pf}

Assume that $d \neq 1,3$. Since the Picard modular group $\conj{\Gam}_3$ has covolume $\frac{1}{216}$, Theorem \ref{volume upper bound} implies that $\Del$ can have smaller covolume than $\conj{\Gam}_3$ if and only if
\[
|d_k| \leq 5.
\]
This excludes all the cases for which Theorem \ref{volume upper bound} holds, leaving $d = 1, 3$. A direct computation shows that no lattice commensurable with $\Gam_1$ can have smaller covolume than $\Gam_3$. The minimal volume lattices for $d = 1$ are generated by the center and either of $\Gam_1$ or its sister. It remains to show that no lattice commensurable with $\Gam_3$ can have smaller covolume than $\conj{\Gam}_3$.

We return to the above estimate:
\[
\mu(\conj{\Gam} \bs \SU(2,1)) \geq -\frac{\prod_{p \in \mc{P}} \lam_p}{3 h_{k,3}} \frac{\mathrm{L}_{\Q(\sqrt{-3})}(-2)}{48}.
\]
Again, $\lam_p / 3 > 1$ for all $p$ at which $\Gam_{K_f}$ is Iwahori, so we can assume it is always of the form $K_p^{v_j}$. For all other $p$, $\lam_p / 3 \geq 1$ with equality if and only if $p = 2$. Therefore, the only cases to consider are where $K_2 = K_p^{v_1}$ and $K_2, K_3 = K_p^{v_1}$. There, Proposition \ref{normalizer index} implies that
\[
\mu(\conj{\Gam} \bs \SU(2,1)) \geq \mu(\Gam_3 \bs \SU(2,1)) > \mu(\conj{\Gam}_3 \bs \SU(2,1)).
\]
This proves the theorem.

\end{pf}



\subsection{}\label{orbifolds}


Consider $\Gam_d = \Gam_{\mathrm{std}, d}$ as in $\S$\ref{minimal} and the surface $\Gam_d \bs \mathbf{H}_\C^2$. Its orbifold Euler-Poincar\'{e} characteristic is $3 \mu(\Gam_d \bs \SU(2,1)) = -\mathrm{L}_k(-2) / 16$ by Theorem \ref{prasad's formula}. The Chern--Gauss--Bonnet formula then implies that $\mathrm{Vol}(M) = 8 \pi^2 \chi(M) / 3$, which gives the following.


\begin{cor}\label{minimal orbifolds}

There are exactly two noncompact arithmetic complex hyperbolic $2$-orbifolds of minimal volume $\pi^2 / 27$. Their fundamental groups are $\conj{\Gam}_3 = \PU(2,1;\mc{O}_3)$ and its sister.

\end{cor}



\begin{rem}

The normalizer of the sister to $\conj{\Gam}_3$ is described explicitly in \cite{Parker}. This is the group that Parker calls $G_2$. Parker's description is morally the same as the one given in $\S$\ref{arithmetic} above.

\end{rem}



\section{Minimal volume manifolds}\label{manifolds}



\subsection{}\label{Gamma(3)}


Let $M$ be a complex hyperbolic manifold. If $M$ is compact, then $\chi(M) \geq 3$. The compact complex hyperbolic manifolds of Euler characteristic three have been completely classified by Prasad--Yeung \cite{Prasad--Yeung} and Cartwright--Steger \cite{Cartwright--Steger}. If $M$ has cusps, then $\chi(M) \geq 1$. For this section, let $\Gam_d$ be the Picard modular group $\PSU(h, \mc{O}_d)$ for
\[
h = \begin{mat}
0 &
0 &
1 \\

0 &
1 &
0 \\

1 &
0 &
0
\end{mat}
\]
and $\conj{\Gam}_d$ its normalizer in $\PU(2,1)$. If $k = \Q(\sqrt{-d})$ has class number one, then $\conj{\Gam}_d = \PU(h, \mc{O}_d)$ and $\Gam_d = \conj{\Gam}_d$ except when $d = 3$. Let $O_d$ and $\wt{O}_d$ denote the respective complex hyperbolic orbifolds.


\begin{ex}\label{parkers manifold}

John Parker \cite{Parker} noticed that $\conj{\Gam}_3$ contains a torsion-free subgroup of index 72. Thus it has Euler-Poincar\'{e} characteristic one and volume $8 \pi^2 / 3$.

\end{ex}



\begin{thm}\label{minimal manifolds}

Let $M$ be a complex hyperbolic $2$-manifold of Euler characteristic one. Then $M$ covers $\wt{O}_3$ or its sister orbifold $\wt{O}_3^{sis}$.

\end{thm}


The proof of Theorem \ref{minimal manifolds} occupies the remainder of this section. Appendix I lists all the manifolds of Euler characteristic one that cover both $\wt{O}_3$ and $\wt{O}_3^{sis}$. That is, their fundamental groups are all conjugate into $\conj{\Gam}_3 \cap \conj{\Gam}_3^{sis}$.


\subsection{}\label{prasad and torsion}



\begin{prop}\label{prasad for euler characteristic one}

Let $M$ be an arithmetically defined cusped complex hyperbolic manifold of Euler characteristic one, and let $k = \Q(\sqrt{-d})$ be the imaginary quadratic field over which $M$ is defined. Then $d = 1$ or $d = 3$.

\end{prop}



\begin{pf}

Proposition 3.5 of \cite{Prasad--Yeung} gives a list
\[
d = 1, 2, 3, 5, 6, 7, 11, 15, 19, 23, 31
\]
of $d$ for which the Euler--Poincar\'{e} characteristic of $\Gam \bs \SU(2,1)$ can be one. Note that $h_{k,3} = 1$ for all but $d = 23,31$.

Let $\Gam$ be the fundamental group of an arithmetic complex hyperbolic manifold of Euler characteristic one and $\conj{\Gam}$ the maximal lattice in $\PU(2,1)$ containing it. Then $\chi(\conj{\Gam} \bs \mathbf{H}_\C^2)$ is $1 / n$ for some $n$ (cf.\ \cite{Prasad--Yeung}).

Let $O_d$ be the Picard modular surface $\Gam_d \bs \mathbf{H}_\C^2$. Note that every $\Gam_d$ contains torsion. In particular, each contains the order two element $m_2$ which lifts to
\[
M_2 = \begin{mat}
-1 &
0 &
0 \\

0 &
1 &
0 \\

0 &
0 &
-1
\end{mat}.
\]
We now eliminate each $d \neq 3,4$ on a case-by-case basis.

\begin{itemize}

\item[$\bullet$] $\mathbf{d = 2}$. The Picard modular surface $O_2$ and its sister are minimal and have Euler--Poincar\'{e} characteristic $\chi(O_2) = 3 / 16$, so they cannot cover a surface of Euler characteristic one.

For any other $\Gam_{K_f}$,
\[
\frac{\prod \lam_p}{16} = \frac{1}{n}
\]
for some $n$. However, no $\lam_p$ divides $16$, so no other maximal lattice can cover an orbifold of Euler--Poincar\'{e} characteristic one.

\item[$\bullet$] $\mathbf{d = 5}$. The surface $O_5$ and its sisters are minimal and have Euler--Poincar\'{e} characteristic $15 / 8$. Every other minimal orbifold has $\chi > 1$, since each $\lam_p$ is greater than three.

\item[$\bullet$] $\mathbf{d = 6}$. Here, $O_6$ and its sisters have $\chi = 23 / 8$, and they are minimal. Since each $\lam_p$ is greater than three, $d = 6$ is impossible.

\item[$\bullet$] $\mathbf{d = 7}$. The surface $O_7$ and its sister have $\chi = 3 / 7$. They are minimal, and cannot cover any space of Euler characteristic one.

Also, $\lam_p = 7$ for $K_3 = K_3^{v_1}$ or $K_3, K_7 = K_p^{v_1}$, and the normalizer is potentially a degree three extension, giving the possibility of $\chi = 1$. The presence of $m_2$ excludes these cases. Since $\lam_p > 7$ for all other $p$, this rules out $d = 7$.

\item[$\bullet$] $\mathbf{d = 11}$. Since $O_{11}$ and its sister have $\chi = 3 /8$, it remains to consider whether it is possible to have
\[
\frac{\prod \lam_p}{8} = \frac{1}{n}.
\]
However, there are no $p$ for which $\lam_p = 2,4,8$.

\item[$\bullet$] $\mathbf{d = 15}$. The surfaces $O_{15}$ and its sisters have Euler--Poincar\'{e} characteristic one, but contain torsion. There is also the possibility that $\Gam_{K_f}$ has index three in its normalizer for $K_2 = K_2^{v_1}$ or its sisters, but each lattice contains torsion.

\item[$\bullet$] $\mathbf{d = 19}$. Since $O_{19}$ and its sister have $\chi = 11/8$ and $\lam_p \geq 3$ for all $p$, this case is impossible.

\item[$\bullet$] $\mathbf{d = 23}$. Here $\chi = 3$ for $O_{23}$ and its sister, but the corresponding lattices have index $3$ in their normalizers, since $\Q(\sqrt{-23})$ has class number $3$. Since $m_2$ is an element of each of these lattices, they are not torsion-free. They are also the unique minimal elements of the commensurability class.

\item[$\bullet$] $\mathbf{d = 31}$. For $O_{31}$ and its sister, $\chi = 6$. Even though the class number of $\Q(\sqrt{-31})$ is three, torsion does not allow for a manifold of Euler characteristic one.

\end{itemize}

\end{pf}



\subsection{}\label{d = 1 computations}


Let $\Gam_1$ denote the Picard modular group $\PU(2, 1; \Z[i])$ for $h$ as in $\S$\ref{Gamma(3)}. It is maximal in $\PU(2,1)$. Here we rule out $\Gam_1$ from being commensurable with a complex hyperbolic 2-manifold of Euler characteristic one. Since $\chi(O_1) = 1 / 32$, it first must be shown that $\Lam_1$ and its sister have no torsion-free subgroup of index 32.

Both $\Gam_1$ and $\Gam_1^{sis}$ contain a subgroup isomorphic to $\Z / 3 \Z$. Explicitly, we can take the subgroups with generators
\[
\begin{mat}
0 &
0 &
i \\

0 &
1 &
0 \\

i &
0 &
-1
\end{mat}\quad\textrm{and}\quad\begin{mat}
0 &
0 &
i/2 \\

0 &
1 &
0 \\

2i &
0 &
-1
\end{mat},
\]
respectively. Therefore, the index in $\Gam_1$ or $\Gam_1^{sis}$ of a torsion-free subgroup must be divisible by 3. This excludes either from containing a torsion-free subgroup of index 32.

It remains to consider other maximal lattices commensurable with $\Gam_1$. Suppose some other maximal lattice commensurable with $O_1$ covers a manifold of Euler characteristic one. Then
\[
\frac{\prod \lam_p}{96} = \frac{1}{n}
\]
for some $n$. In particular, $\lam_p$ divides $96$. No $\lam_p$ for $\Q(i)$ divides $96$, so $O_1$ is not commensurable with a cusped complex hyperbolic manifold of volume $8 \pi^2 / 3$.


\subsection{}\label{d = 3 computations}


Now, consider $\Gam_3$, and let $O_3$ be the corresponding complex hyperbolic orbifold. Recall that $O_3$ has a three-fold quotient $\wt{O}_3$ with Euler--Poincar\'{e} characteristic $1 / 72$ and fundamental group $\PU(2,1;\mc{O}_3) = \conj{\Gam}_3$. If any other maximal orbifold $O \neq \wt{O}_3$ covers a manifold of Euler characteristic one, then
\[
\frac{\prod \lam_p}{72} = \frac{1}{n}
\]
for some integer $n$, which implies that $O$ is either the sister $\wt{O}_3^{sis}$ to $\wt{O}_3$ or is $\conj{\Gam}_{K_f}$, where we have chosen $K_p^{v_1}$ at $\{2\}$ or $\{2, 3\}$.

To rule out $\conj{\Gam}_{K_f}$, where $K_f$ is not hyperspecial at precisely $\{2\}$ or $\{2, 3\}$, let $G$ be the finite subgroup of $\conj{\Gam}_3$ of order 72 \cite[Proposition 3.4.4]{Holzapfel86}. (The order of this group is given as 48 in \cite{Holzapfel86}, but it is a $\Z/2\Z$ extension of $(\Z/6\Z)^2$, thus has order 72.) This determines a subgroup of
\[
\conj{\Gam}_{K_f} \cap \conj{\Gam}_3
\]
of order 36, since the scalar action of the order six elements preserves the corresponding $\mc{O}_k$-lattices. The remaining order two element does not preserve these lattices. This implies that a torsion-free subgroup of $\conj{\Gam}_{K_f}$ must have index at least 36. Therefore the corresponding manifold must have Euler characteristic at least 3. This rules out these two $K_f$s, and completes the proof of Theorem \ref{minimal manifolds}. \qed


\begin{rem}

We close with a remark about nonarithmetic lattices. Deligne--Mostow (cf.\ \cite{Mostow}) and, independently, Thurston \cite{Thurston}, constructed several nonarithmetic lattices in $\SU(2,1)$. Two of them, $\mu_{71}$ and $\mu_{73}$ in \cite{Thurston}, are nonuniform. Via their descriptions as partially compactified moduli of 5 points on the sphere, one can directly compute their Euler characteristics, which are 109/96 and 227/144, respectively. Therefore, their volumes are greater than the volume of the smallest arithmetic orbifolds, and they cannot possibly cover a manifold of Euler characteristic one. Julien Paupert informed me that he and John Parker have constructed new noncompact, nonarithmetic complex hyperbolic 2-orbifolds.

\end{rem}



\subsection*{Acknowledgments}

I first want to thank my advisor, Alan Reid, for countless conversations and his encouragement. I also want to thank Gopal Prasad, Misha Belolipetsky, Daniel Allcock, Julien Paupert, and Vincent Emery for conversations related to this paper and Ted Chinburg and the University of Pennsylvania for their hospitality while I made final revisions. Finally, I thank the referee for very helpful comments and corrections.



\section*{Appendix: Eight cusped complex hyperbolic 2-manifolds with Euler characteristic one}


As in \cite{Conder--Maclachlan}, one can use MAGMA \cite{Bosma--Cannon--Playoust} to find subgroups of small index that do not contain conjugates of some finite collection of elements. Here, the elements we exclude are the conjugacy classes of torsion in $\conj{\Gam}_3$, which are listed in \cite[Proposition 3.4.4]{Holzapfel86}. The presentation used is the two-generator presentation of Falbel--Parker \cite{Falbel--Parker}. We compute the index 4 subgroup which is the intersection of the two maximal lattices, then all index 18 subgroups of it.

We list generators of each lattice as words in Falbel--Parker's generators $J$ and $R_1$ and the integral homology of the manifold. As some relations are quite long, a MAGMA file is available on the author's website which generates these groups, along with presentations and sufficient data to determine the number of cusps. We remark that the two pairs of manifolds with the same $1^{st}$ homology  can be distinguished by computing the number of subgroups of index four in the lattices.


\begin{center}
\begin{longtable}{|c|c|}
\hline
 & $R_1^{-2} J R_1^2 J^{-1}$ \\
Generators & $J R_1^{-1} J R_1^2 J^{-1} R_1^{-1} J$ \\
 & $R_1^2 J^{-1} R_1 J R_1^2 J^{-1} R_1$ \\
\hline
$1^{st}$ homology & $(\Z / 3 \Z)^2 \oplus (\Z / 9 \Z)$ \\
\hline
\# cusps & 4 \\
\hline \hline
 & $R_1^{-1} J R_1 J^{-1} R_1^{-1} J^{-1}$, $J R_1^{-1} J R_1 J^{-1} R_1^{-1} J$ \\
Generators & $R_1^2 J R_1 J R_1^{-1} J^{-1} R_1^{-1}$ \\
 & $R_1 J^{-1} R_1 J R_1^{-3} J^{-1}$, $R_1 J^{-1} R_1^{-1} J R_1 J$ \\
\hline
$1^{st}$ homology & $(\Z / 2 \Z)^4 \oplus (\Z / 4 \Z)$ \\
\hline
\# cusps & 4 \\
\hline \hline
 & $J R_1 J R_1^2 J^{-1} R_1$ \\
Generators & $R_1 J R_1 J R_1^{-1} J R_1^{-1}$, $J^{-1} R_1 J R_1^2 J^{-1} R_1 J^{-1}$ \\
 & $J^{-1} R_1 J^{-1} R_1^2 J R_1 J^{-1}$ \\
\hline
$1^{st}$ homology & $\Z^4$ \\
\hline
\# cusps & 4 \\
\hline \hline
 & $R_1^2 J R_1^{-2} J^{-1}$ \\
Generators & $J^{-1} R_1 J R_1^2 J^{-1} R_1^{-1}$ \\
 & $R_1 J^{-1} R_1^{-1} J^{-1} R_1 J^{-1} R_1^{-1} J R_1^{-2}$ \\
\hline
$1^{st}$ homology & $(\Z / 3 \Z) \oplus \Z^2$ \\
\hline
\# cusps & 2 \\
\hline \hline
 & $J R_1^2 J$ \\
Generators & $J^{-1} R_1 J^{-1} R_1^2 J R_1 J^{-1}$ \\
 & $J R_1^{-2} J R_1 J R_1^2 J^{-1} R_1$ \\
\hline
$1^{st}$ homology & $\Z^2$ \\
\hline
\# cusps & 2 \\
\hline
\end{longtable}
\begin{longtable}{|c|c|}
\hline
 & $R_1^2 J R_1^2$ \\
Generators &  \\
 & $J R_1 J^{-1} R_1^2 J R_1 J^{-1}$ \\
\hline
$1^{st}$ homology & $(\Z / 3 \Z) \oplus (\Z / 9 \Z)$ \\
\hline
\# cusps & 2 \\
\hline \hline
 & $J R_1 J R_1^2 J^{-1} R_1^{-1}$ \\
Generators & $R_1 J R_1 J R_1^{-1} J R_1$ \\
 & $R_1 J R_1 J^{-1} R_1^2 J^{-1} R_1^{-1} J R_1^{-1}$ \\
\hline
$1^{st}$ homology & $(\Z / 3 \Z) \oplus \Z^2$ \\
\hline
\# cusps & 2 \\
\hline \hline
 & $R_1 J R_1^2 J^{-1} R_1 J$ \\
Generators & $R_1^{-1} J R_1 J R_1^{-1} J R_1^{-1}$ \\
 & $R_1 J^{-1} R_1 J R_1^2 J R_1^{-1} J^{-1} R_1^{-1}$ \\
\hline
$1^{st}$ homology & $\Z^2$ \\
\hline
\# cusps & 2 \\
\hline
\end{longtable}
\end{center}



\bibliography{picardvolumes1}


\end{document}